\documentclass[trsc,nonblindrev]{informs3} 
\usepackage{nomencl}
\usepackage{booktabs}
\usepackage{colortbl}
\usepackage{multirow}
\usepackage{lscape}
\usepackage{etoolbox}
\OneAndAHalfSpacedXI 

\usepackage{natbib}
 \bibpunct[, ]{(}{)}{,}{a}{}{,}%
 \def\bibfont{\small}%
 \def\BIBand{and}%

\TheoremsNumberedThrough     

\EquationsNumberedThrough    

\begin{document}

\setcounter{page}{1}
\VOLUME{00}%
\NO{0}%
\MONTH{Xxxxx}
\YEAR{0000}
\FIRSTPAGE{000}%
\LASTPAGE{000}%
\SHORTYEAR{00}
\ISSUE{0000} %
\LONGFIRSTPAGE{0001} %
\DOI{10.1287/xxxx.0000.0000}%

\RUNAUTHOR{Zhu et al.}

\RUNTITLE{A Strong Formulation for Stochastic Multiple Constrained Resources Air Traffic Flow Management with Reroutes}

\TITLE{A Strong Formulation for Stochastic Multiple Constrained Resources Air Traffic Flow Management with Reroutes}

\ARTICLEAUTHORS{%
\AUTHOR{Guodong Zhu, Peng Wei}
\AFF{Department of Aerospace Engineering, Iowa State University, \{\EMAIL{gzhu@iastate.edu},  \EMAIL{pwei@iastate.edu}\}}
} 

\ABSTRACT{%
This paper addresses the air traffic flow management research problem of determining reroute, ground delay and air delay for flights using stochastic weather forecast information. The overall goal is to minimize system-wide reroute and delay costs. This problem is a primary concern in United States and especially in its northeastern region, and is also the key in enhancing the performance of the new FAA Traffic Management Initiative called Collaborative Trajectory Options Program (CTOP). In this work we present six stochastic integer programming models, including two two-stage models and four multistage models, which are based on a notable strong deterministic flight-by-flight level traffic flow formulation and an aggregate Eulerian traffic flow formulation. Our preliminary numerical results show that completely integer solutions can be achieved from linear relaxation, for all proposed models and for both no-route and reroute cases.

}%


\HISTORY{}

\maketitle

%


\makeatletter
\newcommand*{\rom}[1]{\expandafter\@slowromancap\romannumeral #1@}
\makeatother

\renewcommand\nomgroup[1]{%
  \item[\bfseries
  \ifstrequal{#1}{A}{Input Parameters}{%
  \ifstrequal{#1}{B}{Primary Decision Variables}{%
  \ifstrequal{#1}{C}{Auxiliary Variables}
  }}%
]}

\makenomenclature
\mbox{}

\nomenclature[A0]{$\mathcal{R}$}{Set of resources, including departure airports and PCAs} 
\nomenclature[A1]{$\mathcal{C}$}{Set of ordered pairs of resources. $(r,r^\prime)\in \mathcal{C}$ iff $r$ is connected to $r^\prime$ in the directed graph of departure airports and PCAs}
\nomenclature[A2]{$\Delta^{r,r^\prime}$}{ Number of time periods to travel from resource $r$ to $r^\prime$, defined for all pairs $(r,r^\prime) \in \mathcal{C}$ }
\nomenclature[A3]{$Q$}{Set of scenarios, $q=1,\cdots, "|Q"|$}
\nomenclature[A4]{$F$}{Set of flights, $i=1,\cdots,"|F"|$}
\nomenclature[A5]{$p_q$}{Probability that scenario $q$ occurs}
\nomenclature[A6]{$M_{t,q}^r$}{Real capacity of PCA $r$ in time period $t$ under scenario $q$}
\nomenclature[A7]{$T^r_{ij}$}{Set of allowed time periods for flight $i$ taking route $j$ to departs from/pass through airport/PCA $r$}
\nomenclature[A8]{$\underline{T}^r_{ij}$}{First time period in the set $T^r_{ij}$}
\nomenclature[A9]{$\overline{T}^r_{ij}$}{Last time period in the set $T^r_{ij}$}
\nomenclature[Aa]{$\text{Dep}_i$}{Original scheduled departure time of flight $i$ }
\nomenclature[Ab]{$t^r_{ij}$}{Time period in which flight $i$ taking route $j$ is scheduled to cross PCA $r$}
\nomenclature[Ac]{$\Omega_{ij}$}{Ordered set of indices of the airport/PCAs which flight $i$ passes if taking route $j$}
\nomenclature[Ad]{$\Omega^k_{ij}$}{The $k$-th resource along route $j$ of flight $i$}
\nomenclature[Ae]{$N_{ij}$}{Number of PCAs along route $j$ of flight $i$}
\nomenclature[Af]{$c_{ij}, c_g, c_a$}{Cost for flight $i$ taking route $j$, cost for unit ground delay and air delay}
\nomenclature[Ag]{$B$}{Set of branches in the scenario tree, $b = 1,\dots, "|B"|$}
\nomenclature[Ah]{$N_b$}{Number of scenarios corresponding to branch $b$}
\nomenclature[Ai]{$o_b, \mu_b$}{Start and end nodes of branch $b$ }
\nomenclature[Aj]{$\Phi_k$}{Set of indices of the routes which cross PCA $k$}
\nomenclature[Ak]{$\mathcal{P}$}{Set of paths, $\rho=1,\cdots, "|\mathcal{P}"|$}
\nomenclature[Al]{$\rho_1,\rho_{-1}$}{First and last PCA on path $\rho$}

\nomenclature[B0]{$w^{r}_{ijt}$}{Whether flight $i$ taking route $j$ departs from/passes through airport/PCA $r$ by time $t$}
\nomenclature[B1]{$w^{rq}_{ijt}$}{Whether flight $i$ taking route $j$ departs from/passes through airport/PCA $r$ by time $t$ under scenario $q$}
\nomenclature[B2]{$\tilde w^{rq}_{ijt}$}{Whether flight $i$ taking route $j$ reaches the first PCA $r$ on route $j$ by time $t$ under scenario $q$}

\nomenclature[C0]{$\delta_{ij}$}{Binary indicator whether flight $i$ will take route $j$}
\nomenclature[C1]{$\delta_{qtij}$}{Binary indicator whether flight $i$ will depart in time period $t$ and take route $j$ under scenario $q$}
\nomenclature[C2]{$\tilde \delta_{qij}$}{Binary indicator whether flight $i$ will take route $j$ under scenario $q$}
\nomenclature[C3]{$P^k_{t,\rho}$}{ Planned direct demand at PCA $k$ in time period $t$ from flights with same path $\rho$ under scenario $q$}
\nomenclature[C4]{$P^{k,q}_{t,\rho}$}{ Planned direct demand at PCA $k$ in time period $t$ from flights with same path $\rho$ under scenario $q$ }
\nomenclature[C5]{$L^k_{t,\rho,q}$ }{ Number of flights with same path $\rho$  that actually cross PCA $k$ in time period $t$ under scenario $q$ }
\nomenclature[C6]{$A^{k,q}_{t,\rho}$}{ Number of flights with same path $\rho$ taking air delay before entering PCA $k$ }

\printnomenclature

\section{Introduction}
The goal of air traffic flow management is to alleviate projected demand-capacity imbalances at airports and in en route airspace. As a new tool in the Federal Aviation Administration (FAA) NextGen portfolio, Collaborative Trajectory Options Program (CTOP) enables air traffic managers to control traffic through multiple congested airspace regions with a single program, which allows traffic to be managed in an integrated and coordinated way. CTOP also allows airline flight operators to submit a set of reroute options (called a Trajectory Options Set or TOS), which provides great flexibility and efficiency to airspace users.

This paper aims to answer the following research question: given reroute options and probabilistic weather forecast information, what is the theoretical best system performance we can achieve in terms of total route and delay costs? This research question is important in designing CTOP program and in analyzing CTOP performance. It is also rather general and fundamental, and can be meaningful for researchers in other countries. In this work, we will tackle this research question in the stochastic programming framework, in which probabilistic weather forecast will be translated into scenario-based capacity data.

There are two paradigms of air traffic flow management models: \emph{Lagrangian models}, which work at flight-by-flight level and provide trajectories and departure times for each flight, and \emph{Eulerian models}, which work at the aggregate flow-based level and provide counts of aircraft in airspace regions. The main advantage of Lagrangian models is their high flexibility and the ability to cope with flight specific differences. The main advantage of Eulerian models is that the model depends only on the size of the geographic regions of interest rather than on the number of aircraft in the region, and thus can usually be very efficiently solved. In this work we present six stochastic integer programming models to find the optimal delay and reroute policy, including three Lagrangian models and three Lagrangian-Eulerian models, each class of models with varying degree to which traffic managers can modify or revise flights’ controlled departure times and reroute. All the models are at least partially on a flight-by-flight level, because the route choices and reroute costs are different for each flight. The motivations to also formulate Lagrangian-Eulerian class of models are that this class of models have less decision variables and constraints and can potentially take less computation time, and they have the ability to easily control the length of
the queue before a constrained resource, which air traffic controllers care about.

Flight-by-flight level air traffic management models tend to be $\mathcal{NP}$-hard, even in the deterministic case as shown in \cite{bertsimas1998air}. Considering the uncertainty in capacity can only compound the problem. Thus, having a good formulation is crucial to solve realistic size problem instance in real time. The stochastic models presented in this work are based on a famous strong air traffic flow formulation, proposed in \cite{bertsimas1998air} and \cite{bertsimas2011integer}. The Lagrangian-Eulerian models are also partly based on the highly efficiently aggregate  models proposed in \cite{ZhuAggregateJournal}.

To differentiate this work with a previous work (\cite{ZhuDisaggregateJournal}), we will call models in this work binary models, because all decision variables are binary. We will call models in \cite{ZhuDisaggregateJournal} as integer models, since some of the key decision variables are integers.

\section{Preliminary Concepts}\label{section: Preliminary Concepts}
\subsection{Potential Constrained Area and Capacity Scenarios}
\begin{figure}
  \centering
  \includegraphics[width=0.5\linewidth]{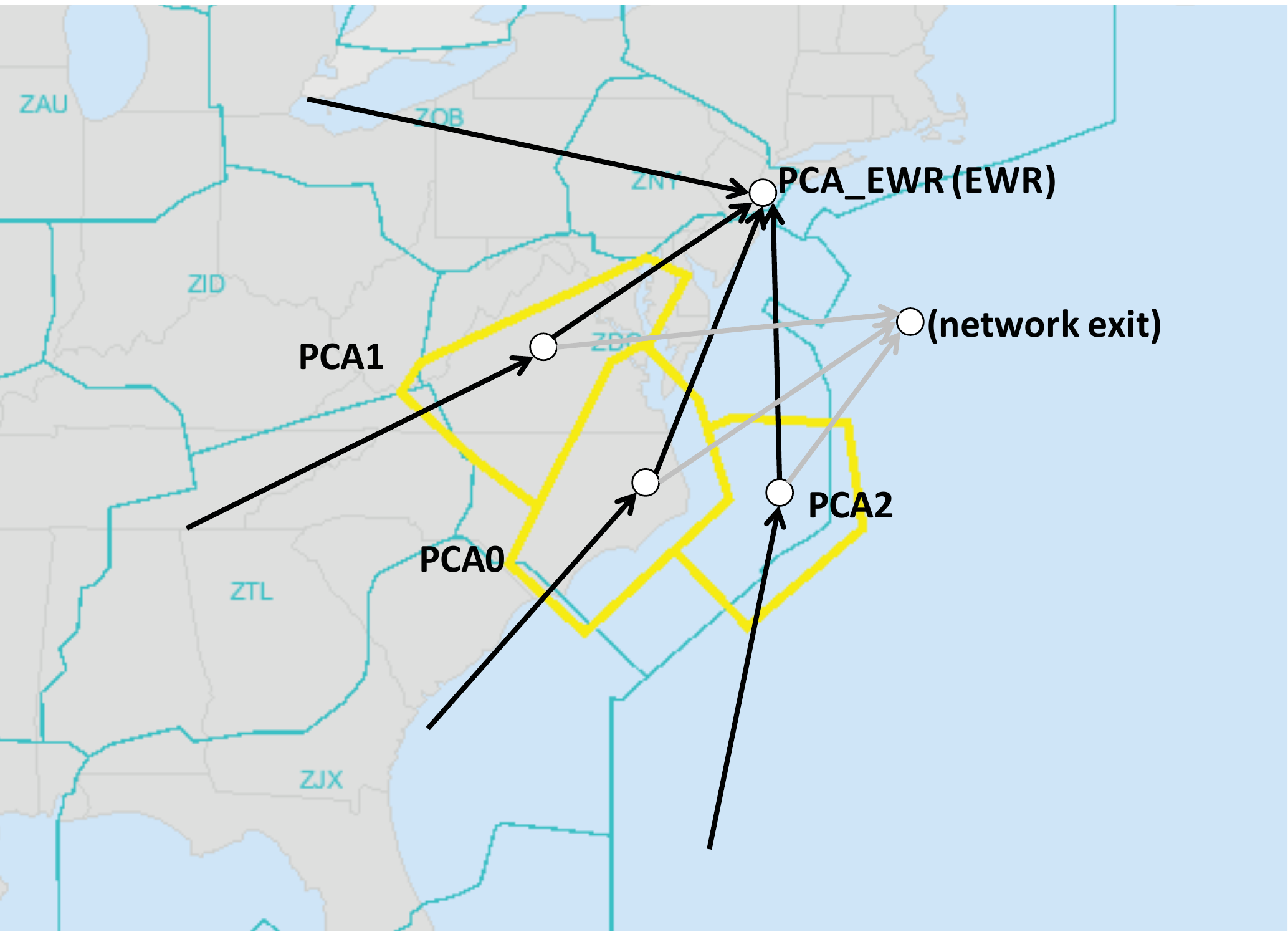}
  \caption{Geographical display of an Airport-PCA Network}\label{PCA_Network}
\end{figure}

\begin{figure}
  \centering
  \includegraphics[width=0.6\linewidth]{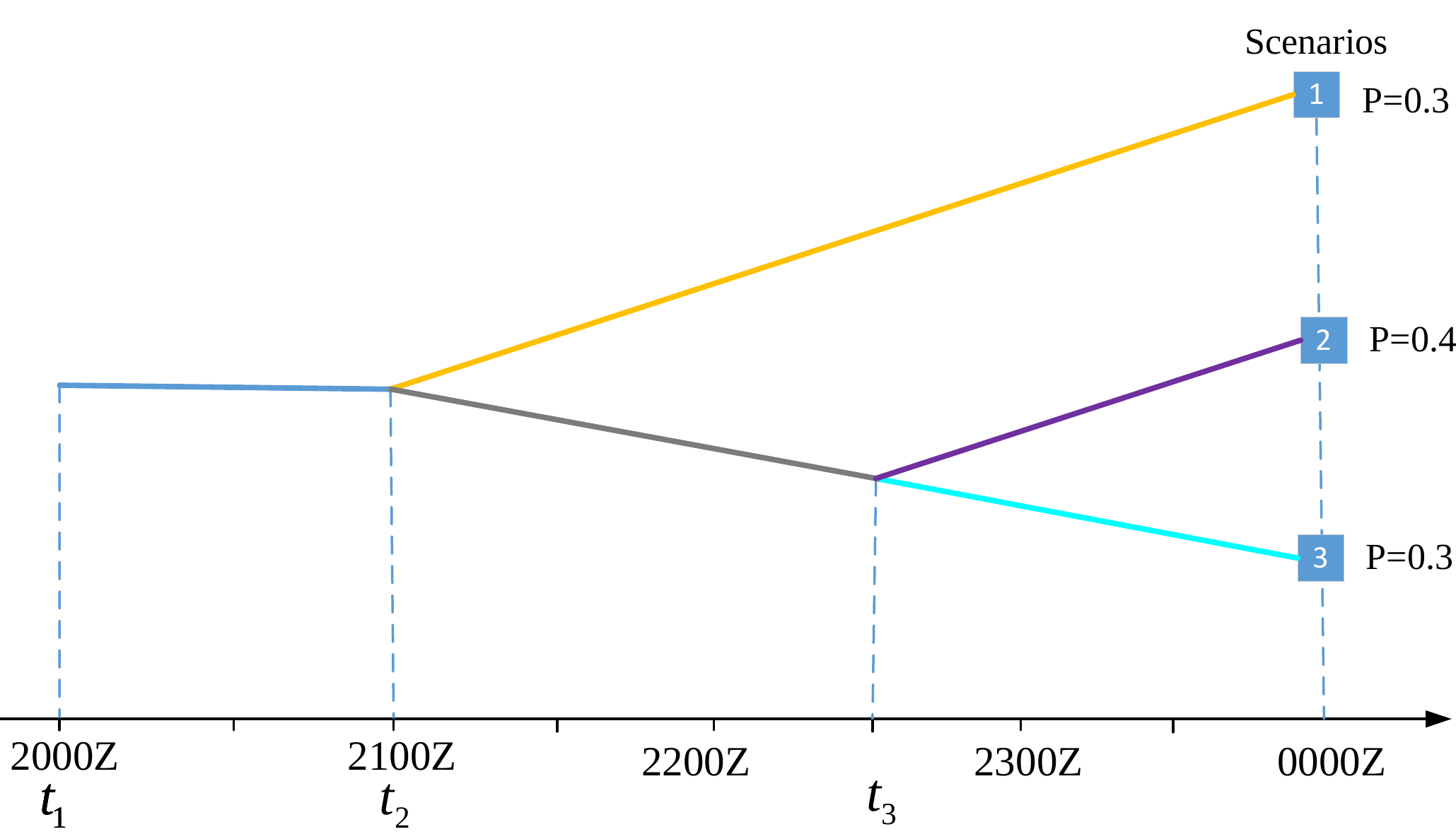}
  \caption{Scenario Tree of PCAs' Evolving Capacities}\label{ScenarioTreeExperiment}
\end{figure}

In this paper, we will model a constrained airspace resource as a Potentially Constrained Area (PCA), in which air traffic demand may exceed capacity and whose future capacity realization is represented by a finite set of scenarios arranged in a scenario tree. A related concept is the airport-PCA network, which refers to a directed graph that links the airports and PCAs, and models the potential movement of traffic between them. Figure \ref{PCA_Network} shows an example of PCA network, which includes three en route PCAs and one constrained airport EWR. Figure \ref{ScenarioTreeExperiment} shows the scenario tree used in this paper. In multi-resource air traffic management problem, the change of operating condition at any PCA will result in a branch point in the scenario tree. Therefore, this scenario tree models the evolution of the future capacities of all four PCAs in Figure \ref{PCA_Network}.

\subsection{Resources along a Route}
In this study, we assume each flight will choose a route from its TOS set, depart from its origin airport, traverse one or more constrained en route PCAs, and land at destination airport (if it is also constrained) or directly exit the PCA system (grey arrow in Figure \ref{PCA_Network}). We define for each flight a list of elements it passes, denoted as $\Omega_{ij}$. By assumption, $\Omega^0_{ij}$ is the departure airport of flight $i$ (for any route $j$) and $\Omega^{k\ge 1}_{ij}$ will include all the PCAs along the route.

\subsection{Path, Direct Demand and Upstream Demand}\label{path}
General multi-resource air traffic management is by nature a multi-commodity problem, since flights will traverse different congested regions and reach different destinations. In our Lagrangian-Eulerian models, after flights choose their routes, they will be grouped by ``path", which is the sequence of PCAs flights traverse. Each path through the PCA nodes in the PCA network establishes a commodity. For example, in Figure \ref{PCA_Network}, $\text{PCA1}\rightarrow \text{PCA}\textunderscore\text{EWR}$ is one path and $\text{PCA1}\rightarrow \text{PCA}\textunderscore\text{Exit}$ is another path. Flights in these two paths share the capacity resource of PCA1. Different routes can correspond to the same path, if they cross the same PCAs.

For each path, we differentiate direct demand, which is the demand for the first PCA on that path and are composed of flights directly flying from departing airports, from upstream demand, which is the demand for the second or later PCAs on that path and are consisting of flights from the upstream en route PCA. We can ground-hold direct demand before flights take off, and we can air-hold both direct demand and upstream demand. Each flight will enter the PCA network through the first PCA (denoted $\rho_1$) on its path $\rho$ and exit through the last PCA (denoted $\rho_{-1}$) on that path. For example, for path $\text{PCA1}\rightarrow \text{PCA}\textunderscore\text{EWR}$, $\rho_1=\text{PCA1}$, $\rho_{-1}=\text{EWR}$.

\section{Two-stage Static Models}\label{Static}
In this section, we introduce the Lagrangian and Lagrangian-Eulerian versions of two-stage stochastic model. In two-stage models, the first stage decisions are the reroute decision and ground delay assignment, and the second stage decisions are the air delays flights need to take in response to the actual weather scenarios.

The primary decision variable in this work is $w^{rq}_{ijt}$, which is a binary variable indicating whether flight $i$ will take $j$ and departs from/passes through airport/PCA $r$ \emph{by time $t$}. To be more clear, when $r$ is an airport ($r=\Omega_{ij}^0$), and if route $j$ is chosen for flight $i$, $w^{rq}_{ijt}=0$ implies that flight is still on the ground. The first time period $w^{rq}_{ijt}=1$ is when this flight is released for departure. When $r$ represents a PCA and $j$ is chosen, $w^{rq}_{ijt}=0$ means flight $i$ is still on its way to PCA $r$, and $w^{rq}_{ijt}$ first becomes 1 when it is admitted to PCA $r$. In two-stage stochastic model, the first stage decisions are made while a flight is still on the ground and are the same for all scenarios, hence we can drop $q$ for $w^{r}_{ij\overline{T}^r_{ij}}$ when $r=\Omega_{ij}^0$.

\subsection{Lagrangian Version}
In the first set of constraints we ensure that one and only route is chosen for each flight:
\begin{equation}\label{only_one_route}
\begin{aligned}
&w^{r}_{ij\overline{T}^r_{ij}} = \delta_{ij} &\quad\forall i\in F, j\in F_i, r=\Omega_{ij}^0\\
&\sum_{j\in F_i}\delta_{ij}=1 &\quad\forall i\in F\\
\end{aligned}
\end{equation}
If $j$ is indeed selected for flight $i$, then this flight must depart by the last allowed departure time period $\overline{T}^r_{ij}$. Here $\delta_{ij}$ is only an ancillary variable.

There are two types of connectivity constraints in this problem: connectivity in time and connectivity between resources. Connectivity between time ensures that if a flight has been admitted to a resource by time $t$, then $w^{r}_{ij,t^\prime}$ has to be 1 for all later time periods $t^\prime > t$.
\begin{equation}
\begin{aligned}
& w^{r}_{ij,t}- w^{r}_{ij,t-1}\ge 0\quad \forall i\in F, j\in F_i, r\in \Omega^0_{ij}, t\in T_{ij}^r, q\in Q \\
& w^{r,q}_{ij,t}- w^{r,q}_{ij,t-1}\ge 0\quad \forall i\in F, j\in F_i, r\in \Omega^{k\ge 1}_{ij}, t\in T_{ij}^r, q\in Q \\
\end{aligned}
\end{equation}
Connectivity between resources impose that if a flight passes through resource $r^\prime$ by $t+\Delta^{r,r^\prime}$, it must have been admitted to $r$, which is the upstream resource on route $j$, by $t$.
\begin{equation}
\begin{aligned}
& w^{r^\prime,q}_{ij,t+\Delta^{r,r^\prime}} - w^{r}_{ij,t} \le 0\quad \forall i\in F, j\in F_i, r=\Omega^{0}_{ij},r^\prime = \Omega^{1}_{ij}, t\in T_{ij}^r\\
& w^{r^\prime,q}_{ij,t+\Delta^{r,r^\prime}} - w^{r,q}_{ij,t} \le 0\quad \forall i\in F, j\in F_i, r,r^\prime \in \Omega^{k\ge 1}_{ij}, t\in T_{ij}^r
\end{aligned}
\end{equation}
The capacity constraint stipulates that the number of flights admitted to PCA $r$ should not exceed its actual capacity at time $t$.
\begin{equation}
\sum_{(i,j)\in \Phi_k; t\in T^r_{ij}} (w^{r,q}_{ijt} - w^{r,q}_{ij, t-1} )\le M_{tq}^r\quad \forall r\in \Omega^{k\ge 1}_{ij}, t\in T, q\in Q
\end{equation}
The boundary conditions are:
\begin{eqnarray}
w^{r}_{ij,\underline{T}^r_{ij}-1} = 0 &\quad \forall i\in F, j\in F_i, r=\Omega_{ij}^0\\
w^{r,q}_{ij, \underline{T}^r_{ij}-1} = 0 &\quad \forall i\in F, j\in F_i, r \in \Omega^{k\ge 1}_{ij}, q\in Q\\
w^{r,q}_{ij, \overline{T}^r_{ij}} = \delta_{ij} &\quad \forall i\in F, j\in F_i, r=\Omega^{N_{ij}}_{ij}, q\in Q
\end{eqnarray}
Ground delay for flight $i$ is:
\begin{equation}
g_{i} = \sum_{j\in F_i}\Big[ \sum_{t\in T^r_{ij}; r = \Omega_{ij}^0}t(w^{r}_{ij,t}- w^{r}_{ij,t-1}) - \delta_{ij}\text{Dep}_i \Big]
\end{equation}
Air delay for flight $i$ under scenario $q$ is:
\begin{equation}
a_{iq} = \sum_{j\in F_i}\Big[ \sum_{t\in T^r_{ij}; r=\Omega_{ij}^{N_{ij}}}\Big(t(w^{r}_{ij,t}- w^{r}_{ij,t-1}) - \delta_{ij}t^{r}_{ij} \Big) \Big] - g_{iq}
\end{equation}
In this work, since we assume flight cannot depart before scheduled departure time and cannot speed up, therefore $\text{Dep}_i=\underline{T}^{\Omega_{ij}^0}_{ij}$, $t^r_{ij}=\underline{T}^{\Omega_{ij}^r}_{ij}$.

The objective function minimizes the total reroute, ground delay, and expected air delay costs. Arranging the terms in the following formula
\begin{equation*}
\min \sum_{i\in F}\Big(c_gg_i+\sum_{q\in Q}c_aa_{iq}+\sum_{j\in F_i}c_{ij}\delta_{ij}\Big)
\end{equation*}
we obtain
\begin{equation}
\begin{split}
& \min\quad \sum_{i\in F}\sum_{j\in F_i}\Big[c_{ij}\delta_{ij}+(c_g-c_a)\sum_{t\in T^r_{ij}; r = \Omega_{ij}^0}\Big( t(w^{r}_{ij,t}- w^{r}_{ij,t-1}) - \delta_{ij}\underline{T}_{ij}^{r}\Big) +  \\
&\quad\quad\quad c_a\sum_{t\in T^r_{ij}; r=\Omega_{ij}^{N_{ij}}}\Big(t(w^{r}_{ij,t}- w^{r}_{ij,t-1}) - \delta_{ij}t^{r}_{ij} \Big) \Big ]\label{TwostageDisaggregateObj}
\end{split}
\end{equation}

\subsection{Lagrangian-Eulerian Version}
The Lagrangian-Eulerian formulation is listed below. As mentioned section \ref{path}, once a flight has chosen a routed, left the airport and \emph{arrives at} the first PCA on the picked route, it will be grouped into traffic flows along that path. That is exactly what constraints (\ref{Disaggregate_Aggregate}) describes. The key word here is \emph{arrives at}, which is different with \emph{pass through}.
\begin{equation}
\begin{split}
\min\quad \sum_{i\in F}\sum_{j\in F_i}\Big[c_{ij}\delta_{ij}+(c_g-c_a)\sum_{t\in T^r_{ij}; r = \Omega_{ij}^0}\Big( t(w^{r}_{ij,t}- w^{r}_{ij,t-1}) - \delta_{ij}\underline{T}_{ij}^{r}\Big) \Big] + c_a\sum_{q\in Q}p_q\sum_{t\in T}\sum_{\rho\in \mathcal{P}}\sum_{k \in \rho}A^{k,q}_{t,\rho}
\end{split}
\end{equation}
\begin{align}
&w^{r}_{ij\overline{T}^r_{ij}} = \delta_{ij} &\quad\forall i\in F, j\in F_i, r=\Omega_{ij}^0\\
&\sum_{j\in F_i}\delta_{ij}=1 &\quad\forall i\in F\\
& w^{r}_{ij,t}- w^{r}_{ij,t-1}\ge 0 &\quad \forall i\in F, j\in F_i, r= \Omega^0_{ij}, t\in T_{ij}^r, q\in Q \\
&\tilde{w}^{r^\prime,q}_{ij,t}- \tilde{w}^{r^\prime,q}_{ij,t-1}\ge 0& \forall i\in F, j\in F_i, r^\prime= \Omega^1_{ij}, t\in T_{ij}^{r^\prime}, q\in Q\\
& w^{r,q}_{ij, \underline{T}^r_{ij}-1} = \tilde{w}^{r^\prime,q}_{ij, \underline{T}^{r^\prime}_{ij}-1} = 0 & \quad \forall i\in F, j\in F_i, r=\Omega^0_{ij}, r^\prime = \Omega^1_{ij}, q\in Q \\
& \tilde{w}^{r^\prime,q}_{ij,t+\Delta^{r,r^\prime}} - w^{r}_{ij,t} = 0  &\quad \forall i\in F, j\in F_i, r=\Omega^{0}_{ij},r^\prime = \Omega^{1}_{ij}, t\in T_{ij}^r\\
& P^k_{t,\rho}= \sum_{(i,j)\in \Phi_k; j\in \rho}\sum_{t\in T^r_{ij}; r = \Omega_{ij}^1}(\tilde{w}^{r}_{ij,t} - \tilde{w}^{r}_{ij,t-1}) &\quad \forall t\in T,\rho\in \mathcal{P}, k=\rho_1 \label{Disaggregate_Aggregate}\\
& L^k_{t,\rho,q}            = \begin{cases}
                             \mbox{if $k=\rho_1$} &   P^k_{t,\rho}  -(A^k_{t,\rho,q}-A^k_{t-1,\rho,q}) \label{Static_L}\\
                             \mbox{else}       &   \text{UpPCA}^k_{t,\rho,q}-(A^k_{t,\rho,q}-A^k_{t-1,\rho,q})\\
                           \end{cases}& \forall t\in T, q\in Q, \rho\in \mathcal{P}, k\in \rho  \\
& \text{UpPCA}^k_{t,\rho,q} = L^{k^\prime}_{t-\Delta^{k^\prime,k},\rho,q} & t\in T, q\in Q, (k^\prime,k)\in \rho \label{Static_UpPCA}\\
& \sum_{t\in T}P^{k=\rho_1}_{t,\rho} = \sum_{t\in T} L^{k=\rho_{-1}}_{t,\rho,q}\quad \forall \rho\in \mathcal{P}, q\in Q \label{bh1}\\
& P^k_{t,\rho}, L^k_{t,\rho,q}, A^k_{t,\rho,q}\ge 0& \forall t\in T, q\in Q, \rho\in \mathcal{P}, k\in \rho \\
& M^k_{t,q}\ge \sum_{\rho\in \mathcal{P}}L^k_{t,\rho,q}& \forall t\in T, q\in Q, k\in P \label{Static_Capacity}
\end{align}

\section{Multistage Dynamic Models}\label{Dynamic}
In this section, we introduce the multistage stochastic models which can dynamically adjust flight release time and reroute choice before actual departure.
\subsection{Lagrangian Version}
The formulation is listed as follows:
\begin{align}
\min\quad & \sum_{q\in Q}p_q\sum_{i\in F}\sum_{j \in F_i}\Bigg[\Big( c_{ij} - (c_g-c_a)\underline{T}_{ij}^{r=\Omega_{ij}^0} - c_a\underline{T}_{ij}^{r=\Omega_{ij}^{N_{ij}}} \Big) \tilde{\delta}_{qij}+ (c_g-c_a) \sum_{t\in T^r_{ij}; r = \Omega_{ij}^0}t(w^{r,q}_{ij,t}- w^{r,q}_{ij,t-1}) + \nonumber \\
&\quad\quad\quad\quad\quad\quad\quad\quad c_a \sum_{t\in T^r_{ij}; r = \Omega_{ij}^{N_{ij}}}t(w^{r,q}_{ij,t}- w^{r,q}_{ij,t-1}) \Bigg]
\end{align}
\begin{align}
&\delta_{qtij} =  w^{r,q}_{ijt} - w^{r,q}_{ij, t-1}       & \forall i\in F, j\in F_i, r = \Omega^0_{ij}, q\in Q\\
&\tilde{\delta}_{qij} = \sum_{t\in T_{ij}^r; r = \Omega^0_{ij}}\delta_{qtij} & \forall i\in F, j\in F_i, q\in Q\\
&\sum_{j\in F_i}\tilde{\delta}_{qij} = 1&\forall i\in F, q\in Q \label{MultistageDelta} \\
&w^{r,q}_{ij,t}- w^{r,q}_{ij,t-1}\ge 0& \forall i\in F, j\in F_i, r\in \Omega_{ij}, t\in T_{ij}^r, q\in Q \label{connectivity_in_time}\\
&w^{r^\prime,q}_{ij,t+\Delta^{r,r^\prime}} - w^{r,q}_{ij,t} \le 0& \forall i\in F, j\in F_i, r,r^\prime \in \Omega_{ij}, t\in T_{ij}^r \label{connectivity_between_resources}\\
&\sum_{(i,j)\in \Phi_k; t\in T^r_{ij}} (w^{r,q}_{ijt} - w^{r,q}_{ij, t-1} )\le M_{tq}^r& \forall r\in \Omega^{k\ge 1}_{ij}, t\in T, q\in Q \label{capacity_constraint}\\
& \delta_{q^b_1tij}=\dots= \delta_{q^b_{N_b}tij}& \forall i\in F,j\in F_i,t \in T^{\Omega_{ij}^0}_{ij}, b\in B,\mu_b\ge t\ge o_b\label{RouteNonanticipativity} \\
& w^{r,q}_{ij, \underline{T}^r_{ij}-1} = 0 & \quad \forall i\in F, j\in F_i, r \in \Omega_{ij}, q\in Q \label{boundary_condition}\\
& w^{r,q}_{ij, \overline{T}^r_{ij}} = \tilde{\delta}_{qij} & \forall i\in F, j\in F_i, q\in Q, r=\Omega^{N_{ij}}_{ij} \label{no_delay_if_not_chosen}
\end{align}
The first three sets of constraints make sure one and only route will be chosen for each flight. $\delta_{qtij}$ is an ancillary binary variable indicating whether flight will take route $j$ and depart in time period $t$. $\tilde \delta_{qij}$ is another ancillary variable which shows whether flight will choose route $j$ under scenario $q$. (\ref{connectivity_in_time}) and (\ref{connectivity_between_resources}) are connectivity in time constraint and connectivity between resources constraint. (\ref{capacity_constraint}) is the capacity constraint, which has exactly the same expression as in the two-stage model. In multistage model, we will also have a set of nonanticipativity constraints (\ref{RouteNonanticipativity}), which ensures that decisions are made solely based on the information available at that time. (\ref{boundary_condition}) and (\ref{no_delay_if_not_chosen}) are boundary conditions.

\subsection{Lagrangian-Eulerian Version}
The dynamic Lagrangian-Eulerian model is straightforward. The "Lagrangian part" is similar to several constraints in dynamic Lagrangian model, and the "Eulerian part" is the exactly same as in the static Lagrangian-Eulerian except for the additional superscript $q$ in $P^{k,q}_{t,\rho}$ and $\tilde{w}^{r^\prime,q}_{ij,t}$.
\begin{equation}
\begin{split}
& \min\quad \sum_{q\in Q}p_q\Bigg[ \sum_{i\in F}\sum_{j\in F_i}\sum_{t \in {T}_{ij}^{r=\Omega_{ij}^0}}\Big(c_g(t-\underline{T}_{ij}^{r=\Omega_{ij}^0})+c_{ij}\Big)\delta_{qtij}     + c_a\sum_{t\in T}\sum_{\rho\in \mathcal{P}}\sum_{k \in \rho}A_{t,\rho}^{k,q} \Bigg]
\end{split}
\end{equation}
\begin{align}
&\delta_{qtij} =  w^{r,q}_{ijt} - w^{r,q}_{ij, t-1}       & \forall i\in F, j\in F_i, r = \Omega^0_{ij}, q\in Q\\
&\tilde{\delta}_{qij} = \sum_{t\in T_{ij}^r; r = \Omega^0_{ij}}\delta_{qtij} & \forall i\in F, j\in F_i, q\in Q\\
&\sum_{j\in F_i}\tilde{\delta}_{qij} = 1&\forall i\in F, q\in Q  \\
&w^{r,q}_{ij,t}- w^{r,q}_{ij,t-1}\ge 0& \forall i\in F, j\in F_i, r= \Omega^0_{ij}, t\in T_{ij}^r, q\in Q\\
&\tilde{w}^{r^\prime,q}_{ij,t}- \tilde{w}^{r^\prime,q}_{ij,t-1}\ge 0& \forall i\in F, j\in F_i, r^\prime= \Omega^1_{ij}, t\in T_{ij}^{r^\prime}, q\in Q\\
&\tilde{w}^{r^\prime,q}_{ij,t+\Delta^{r,r^\prime}} - w^{r,q}_{ij,t} = 0& \forall i\in F, j\in F_i, r=\Omega^{0}_{ij},r^\prime = \Omega^{1}_{ij}, t\in T_{ij}^r\\
& \delta_{q^b_1tij}=\dots= \delta_{q^b_{N_b}tij}& \forall i\in F,j\in F_i,t \in T^{\Omega_{ij}^0}_{ij}, b\in B,\mu_b\ge t\ge o_b \\
& w^{r,q}_{ij, \underline{T}^r_{ij}-1} = \tilde{w}^{r^\prime,q}_{ij, \underline{T}^r_{ij}-1} = 0 & \quad \forall i\in F, j\in F_i, r=\Omega^0_{ij}, r^\prime = \Omega^1_{ij}, q\in Q \\
& \tilde{w}^{r^\prime,q}_{ij, \overline{T}^{r^\prime}_{ij}} = \tilde{\delta}_{qij} & \forall i\in F, j\in F_i, q\in Q, r^\prime=\Omega^{1}_{ij}\\
& P^{k,q}_{t,\rho}= \sum_{(i,j)\in \Phi_k; j\in \rho}\sum_{t\in T^r_{ij}; r^\prime = \Omega_{ij}^1}(\tilde{w}^{r^\prime,q}_{ij,t} - \tilde{w}^{r^\prime,q}_{ij,t-1}) &\quad \forall t\in T,\rho\in \mathcal{P}, k=\rho_1, q \in Q\\
& (\ref{Static_L})-(\ref{Static_Capacity}) \nonumber
\end{align}
\subsection{Multistage Semi-dynamic Models}\label{Semi}
In multi-stage dynamic models, a flight can revise the departure time and reroute choice multiple times as long as it is on the ground. An important and practical variant of multi-stage model is semi-dynamic model, in which the ground delay and reroute decisions are made at some pre-determined time, e.g. 1 hour before schedule departure time. For simplicity, we usually assume the decisions are made the scheduled departure time. Instead of enforcing (\ref{RouteNonanticipativity}), we will impose the following nonanticipativity constraints in Lagrangian and Lagrangian-Eulerian models:
\begin{equation}
\begin{split}
\delta_{q^b_1tij}=\dots= \delta_{q^b_{N_b}tij}&\quad\forall i\in F,j\in F_i, t \in T^{\Omega_{ij}^0}_{ij}, b\in B,\mu_b\ge \text{Dep}_i\ge o_b\\
\end{split}
\end{equation}
The major advantage of semi-dynamic model over dynamic model is the higher predictability in flight schedule.

\section{Experimental Results}\label{Experiment}
\begin{figure}
\centering
\begin{minipage}[t]{.5\textwidth}
  \centering
  \includegraphics[width=\linewidth]{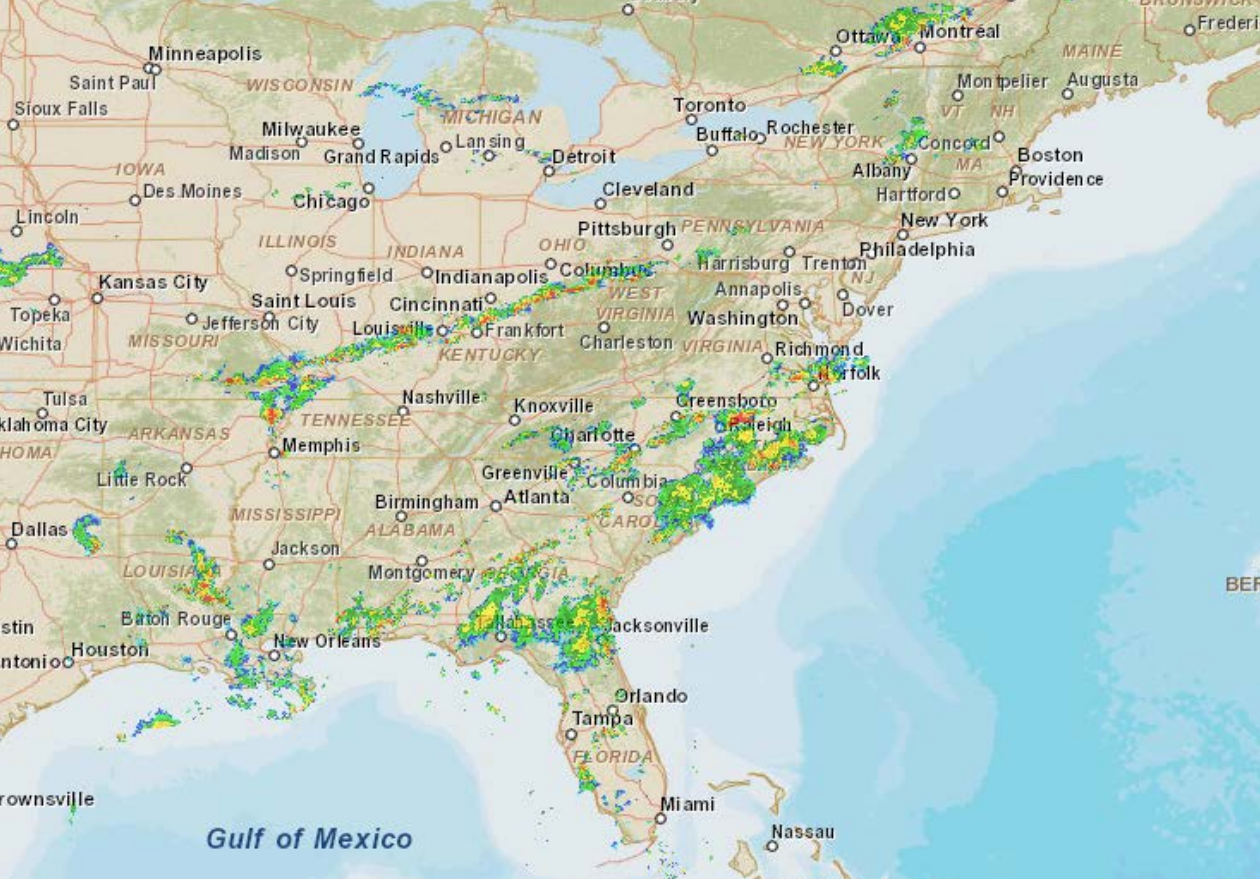}
  \caption{Weather Forecast for 2210z, Taken at 1522z}\label{Convective_Weather_Forecast}
\end{minipage}%
\begin{minipage}[t]{.5\textwidth}
  \centering
  \includegraphics[width=0.85\linewidth]{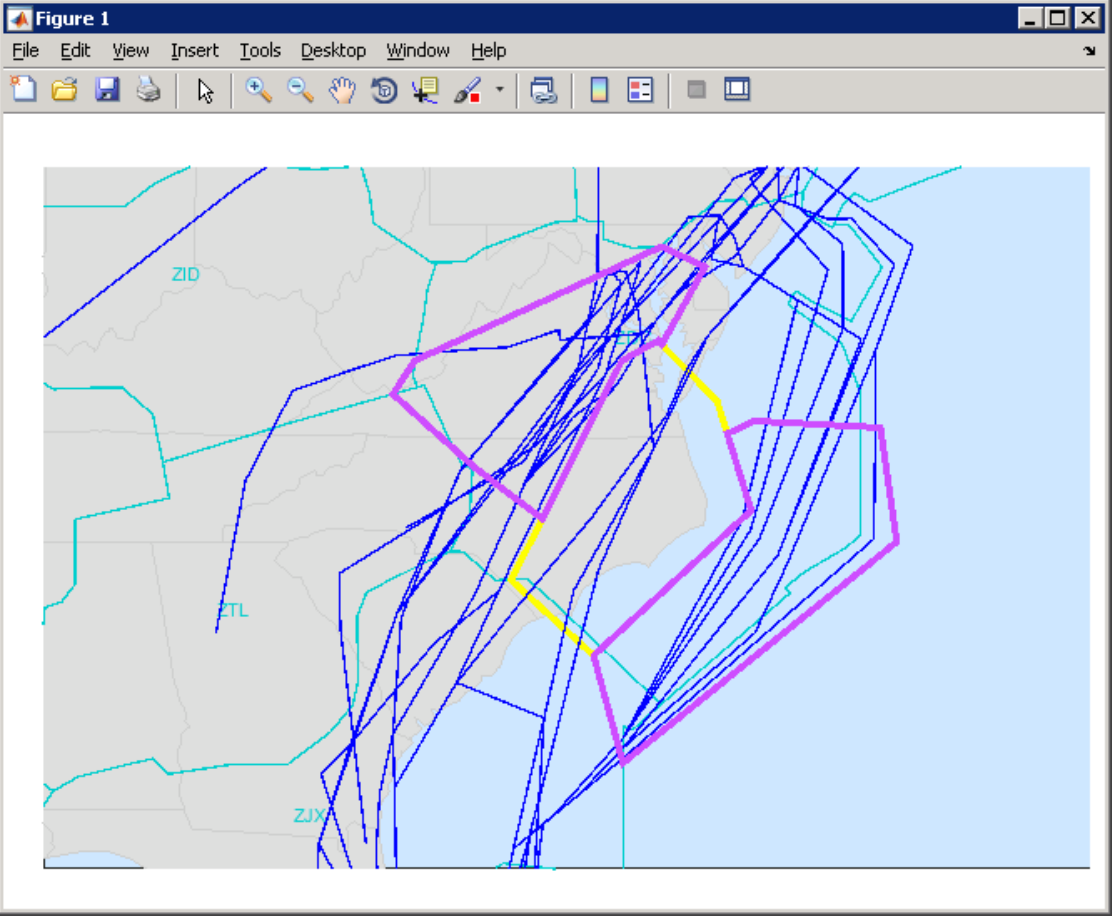}
  \caption{Traffic Routing Around the Original PCA}\label{Traffic_Routing}
\end{minipage}
\end{figure}

To demonstrate the performance of the proposed models, we created an operational use case based on actual events from July 15, 2016. This use case primarily addresses convective weather activity in southern Washington Center (ZDC) and EWR airport. Figure \ref{Convective_Weather_Forecast} shows the pattern of convective weather activity for that day. There is a four-hour capacity reduction in ZDC/EWR from 2000z to 2359z. By analyzing the traffic trajectory (Figure \ref{Traffic_Routing}) and weather data,  we can build the airport-PCA network, shown in Figure \ref{PCA_Network}.
%

\subsection{Capacity Profiles and Traffic Demand}
For comparison purposes, we use the same capacity data as in \cite{ZhuDisaggregateJournal}. The detailed capacity information is listed in Table \ref{Capacity_Scenarios}. The three scenarios correspond to optimistic, average, and pessimistic weather forecast. We can see that in scenario 1 at 2100Z PCA1's 15-minute capacity changes from 44 to 50, the EWR's capacity changes from 8 to 10; in scenario 2 at 2230Z, the capacities of PCA1 and EWR return to the nominal values.  These two changes correspond to the two branch points in the scenario tree shown in Figure \ref{ScenarioTreeExperiment}.

In GDP optimization, we usually add one extra time period to make sure all flights will land at the end of the planning horizon. Because CTOP has multiple constrained resources, we need to add more than one time period depending on the topology of the FCA-PCA network. In this case, we add eight extra time periods, because the longest travel time between the three en route PCAs and EWR among all TOS options is around 2 hours (8 time periods). For any time periods outside the CTOP start-end time, e.g. the eight extra time periods in Table \ref{Capacity_Scenarios}, nominal capacities are used.

We use flight trajectory data from System Wide Information Management (SWIM) and  Coded Departure Route (CDR) database for traffic demand modeling. In total 1098 flights are captured by this CTOP, among them 890 flights that traverse the PCAs in their active periods. And there are in total 1368 TOS options for 890 flights, on average 1.54 options per flight.

\begin{table*}[t]
  \resizebox{1.0\textwidth}{!}{%
  \centering
  \begin{tabular}{cc >{\columncolor[gray]{0.8}}c >{\columncolor[gray]{0.8}}c >{\columncolor[gray]{0.8}} c >{\columncolor[gray]{0.8}}c cccc
  >{\columncolor[gray]{0.8}}c >{\columncolor[gray]{0.8}}c >{\columncolor[gray]{0.8}}c >{\columncolor[gray]{0.8}}c
  cccc
  >{\columncolor[gray]{0.8}}c  >{\columncolor[gray]{0.8}}c >{\columncolor[gray]{0.8}} c >{\columncolor[gray]{0.8}}c c ccc }
  \toprule
  & Resource/Time Bin & 20:00 & 15 & 30 & 45 & 21:00 & 15 & 30 & 45 & 22:00 & 15& 30 & 45 & 23:00 & 15 & 30 & 45 & 00:00 & 15 & 30 & 45 &01:00&15&30&45\\
  \midrule
  \multirow{4}{*}{Scen1} & PCA0  & 13 & 13 & 13 & 13 & 25 & 25 & 25 & 25 & 25 & 25 & 25 & 25 & 25 & 25 & 25 & 25 & 25 & 25 & 25 & 25&25& 25 & 25 & 25\\
                         & PCA1 & 44 & 44 & 44 & 44 & 50 & 50 & 50 & 50 & 50 & 50 & 50 & 50 & 50 & 50 & 50 & 50 & 50 & 50 & 50 & 50 &50& 50 & 50 & 50\\
                         & PCA2 &  5 &  5 &  5 &  5 &  5 &  5 &  5 &  5 &  5 &  5 &  5 &  5 &  5 &  5 &  5 &  5 &  5 &  5 &  5 &  5 & 5&  5 &  5 &  5\\
                         & EWR  &  8 &  8 &  8 &  8 & 10 & 10 & 10 & 10 & 10 & 10 & 10 & 10 & 10 & 10 & 10 & 10 & 10 & 10 & 10 & 10  & 10 & 10 & 10 & 10\\
  \midrule
  \multirow{4}{*}{Scen2} & PCA0  & 13 & 13 & 13 & 13 & 13 & 13 & 13 & 13 & 13 & 13 & 25 & 25 & 25 & 25 & 25 & 25	& 25 & 25 & 25 & 25 & 25& 25 & 25 & 25\\
                         & PCA1 & 44 & 44 & 44 & 44 & 44 & 44 & 44 & 44 & 44 & 44 & 50 & 50 & 50 & 50 & 50 & 50 & 50 & 50 & 50 & 50 & 50 & 50 & 50 & 50\\
                         & PCA2 &  5 &  5 &  5 &  5 &  5 &  5 &  5 &  5 &  5 &  5 &  5 &  5 &  5 &  5 &  5 &  5 &  5 &  5 &  5 &  5 & 5&  5 &  5 &  5 \\
                         & EWR  &  8 &  8 &  8 &  8 &  8 &  8 &  8 &  8 &  8 &  8 & 10 & 10 & 10 & 10 & 10 & 10 & 10 & 10 & 10 & 10 & 10& 10 & 10 & 10 \\
  \midrule
  \multirow{4}{*}{Scen3} & PCA0  & 13 & 13 & 13 & 13 & 13 & 13 & 13 & 13 & 13 & 13 & 13 & 13 & 13 & 13 & 13 & 13 & 25 & 25 & 25 & 25 & 25& 25 & 25 & 25\\
                         & PCA1 & 44 & 44 & 44 & 44 & 44 & 44 & 44 & 44 & 44 & 44 & 44 & 44 & 44 & 44 & 44 & 44 & 50 & 50 & 50 & 50 & 50& 50 & 50 & 50\\
                         & PCA2 &  5 &  5 &  5 &  5 &  5 &  5 &  5 &  5 &  5 &  5 &  5 &  5 &  5 &  5 &  5 &  5 &  5 &  5 &  5 &  5 & 5&  5 &  5 &  5\\
                         & EWR  &  8 &  8 &  8 &  8 &  8 &  8 &  8 &  8 &  8 &  8 &  8 &  8 &  8 &  8 &  8 &  8 & 10 & 10 & 10 & 10 & 10& 10 & 10 & 10 \\
  \bottomrule
  \end{tabular}}
  \caption{Capacity Scenarios} \label{Capacity_Scenarios}
\end{table*}


\subsection{Model Comparisons}
\begin{table}[h!]
\centering
\resizebox{1.0\textwidth}{!}{%
  \begin{tabular}{lp{1.5cm}lrrrrrrr}
    \toprule
    &&&\multicolumn{3}{c}{Lagrangian}&\multicolumn{3}{c}{Lagrangian-Eulerian}\\
    \cmidrule(lr){4-6}\cmidrule(lr){7-9}
    &&&Variables&Constraints&Non-zeros&Variables&Constraints&Non-zeros\\
   \multirow{6}{*}{\parbox{1cm}{Integer Model}}&\multirow{3}{*}{\parbox{1cm}{Without TOS}} & Static       &  53,462  &  8,730  &  152,232  & 17,532 & 3,698 &  51,524      \\
                                              && Semi-dynamic &  55,242  &  9,933  &  154,638    & 49,656     &  8,796      &  145,902          \\
                                              && Dynamic      &  97,584 &  20,892 &  260,358    &  91,998     &  19,758    &   247,668         \\
   \cmidrule{2-9}
   &\multirow{3}{*}{\parbox{1cm}{With TOS}} & Static       & 63,801 & 14,969 & 188,688 &27,049 & 5,144 &    77,605               \\
                                           && Semi-dynamic & 69,273 & 20,389 & 204,176  & 78,087& 19,534 &   237,287            \\
                                          & & Dynamic      & 130,074& 30,463 & 410,111 & 137,703& 25,504&  37,2563           \\
  \midrule
  \midrule
   \multirow{6}{*}{\parbox{1cm}{Binary Model}}&\multirow{3}{*}{\parbox{1cm}{Without TOS}} & Static  & 65,652   & 118,368   & 323,674   &  33,426& 50,490 & 130,104       \\
                                              && Semi-dynamic & 142,452   &  217,763  &  565,696    & 142,350     &   215,579     &  559,966          \\
                                              && Dynamic      & 142,452   &  204,660  &  539,490    & 142,350      &  202,476    &   533,760         \\
   \cmidrule{2-9}
   &\multirow{3}{*}{\parbox{1cm}{With TOS}} & Static       & 100,899 & 181,736 & 499,104  & 51,025& 77,550 &  198,859                \\
                                           && Semi-dynamic & 219,411 & 335,868 & 935,352  & 217,839& 330,672 &  920,433             \\
                                          & & Dynamic      & 219,411 & 315,754 & 895,118  & 217,839& 310,546 &  880,181           \\
    \bottomrule
  \end{tabular}
}
\caption{Model Complexity}\label{Model_Complexity}
\end{table}
\begin{table*}[h!]
  \resizebox{1.0\textwidth}{!}{%
  \centering
    \begin{tabular}{c c c c c c c c c c c c c}
    \toprule
    \multirow{3}{*}{\parbox{2.6cm}{Lagrangian vs. Lagrangian-Eulerian}} & \multicolumn{3}{c}{Ground Delay Periods} & \multicolumn{3}{c}{Air Holding Periods} & &\multicolumn{2}{c}{Integer Model}&Binary Model \\
    \cmidrule(lr){9-10}\cmidrule{11-11}
     & \multicolumn{3}{c}{If This Scenario Occurs:} & \multicolumn{3}{c}{If This Scenario Occurs:} & \multirow{2}{*}{ Expected Cost } & \multirow{2}{*}{ Running Time } & \multirow{2}{*}{ Early Stop} & \multirow{2}{*}{ Running Time }\\
     \cmidrule{2-7}
                            &Scen1 & Scen2 & Scen3  &Scen1 & Scen2 & Scen3 &       & mins & at 1min/3mins  &mins\\
    \midrule
     Two-stage Model        & 284  & 284   &  284   &  0   &  1    &  205  &407.8 & {\cellcolor[gray]{0.8}}0.61 &             & {\cellcolor[gray]{0.8}}0.12 (LP Rlx)\\
     Semi-dynamic Model     & 164  & 285   &  417   & 0    &  0    &  73   &332.1 & {\cellcolor[gray]{0.8}}3.85 &  334.6/332.1& {\cellcolor[gray]{0.8}}0.23 (LP Rlx)\\
     Dynamic Model          & 125  & 288   &  477   &  0   &  0    &  13   &303.6 & {\cellcolor[gray]{0.8}}$>10.0$ & 306.1/303.9  & {\cellcolor[gray]{0.8}}0.21 (LP Rlx)\\
     Perfect Information    & 90   & 285   &  489   &  0   &  0    &   0   &287.7\\
     \midrule
     Two-stage Model        & 284  & 284   &  284   &  0   &  0    &  200  &404.0 & {\cellcolor[gray]{0.8}}0.01 (LP Rlx) & & {\cellcolor[gray]{0.8}}$0.04$ (LP Rlx)\\
     Semi-dynamic Model     & 163  & 284   & 417    & 0    & 0     & 69    &329.0 & {\cellcolor[gray]{0.8}}0.40 &&  {\cellcolor[gray]{0.8}}0.20 (LP Rlx)\\
     Dynamic Model          & 127  & 284   & 472    & 0    & 0     & 12    &300.5 & {\cellcolor[gray]{0.8}}7.69& 300.5&{\cellcolor[gray]{0.8}}0.21 (LP Rlx)\\
     Perfect Information    & 93   & 284   & 484    & 0    & 0     & 0     &286.7 & \\
  \bottomrule
  \end{tabular}}
  \caption{ Stochastic Models Comparison (Delay cost ratio $c_a/c_g=2$) {\bf{WITHOUT}} TOS}\label{WithoutTOSTable}
\end{table*}

\begin{table*}[h!]
  \resizebox{1.0\textwidth}{!}{%
  \centering
    \begin{tabular}{c c c c c c c c c c c c c c c c c c}
    \toprule
    \multirow{3}{*}{} & \multicolumn{3}{c}{RTC Costs in Mins}  & \multicolumn{3}{c}{Ground Delay Periods} & \multicolumn{3}{c}{Air Holding Periods} & &\multicolumn{2}{c}{Integer Model}&Binary Model\\
    \cmidrule(lr){12-13}\cmidrule{14-14}
     & \multicolumn{3}{c}{If This Scenario Occurs:} & \multicolumn{3}{c}{If This Scenario Occurs:} & \multicolumn{3}{c}{If This Scenario Occurs:} &  \multirow{2}{*}{ Expected Cost } & \multirow{2}{*}{ Running Time } & \multirow{2}{*}{ Early Stop} &\multirow{2}{*}{ Running Time }\\
      \cmidrule{2-10}
                        &Scen1     & Scen2    & Scen3     &Scen1       & Scen2    & Scen3    &Scen1    & Scen2    & Scen3        &        &mins & at 1min/3mins&mins\\
    \midrule
     Two-stage Model    &  272      &  272     &   272     &    108    &  108 &   108    & 0   & 0    & 38   & 167.07& {\cellcolor[gray]{0.8}}0.20  & &{\cellcolor[gray]{0.8}}0.19 (LP Rlx)\\
     Semi-dynamic Model &  216      &  256     &   268     &    86     &   108&   130    & 0   & 0    & 22   & 154.21& {\cellcolor[gray]{0.8}}0.08  & &{\cellcolor[gray]{0.8}}0.64 (LP Rlx)\\
     Dynamic Model      &  216      &  256     &   268     &    73      &  110&   136    & 0   & 0    &16    & 149.31& {\cellcolor[gray]{0.8}}1.86  &149.31 &{\cellcolor[gray]{0.8}}0.58 (LP Rlx)\\
     Perfect Information&   76      &  222     &   286     &    56      &  109     &   144    & 0       & 0         &      0       & 129.92\\
  \midrule
       Two-stage Model  & 270  & 270      &  270      &   106      &  106     &  106     &  0       & 0        &   29 &159.40 &{\cellcolor[gray]{0.8}}$<0.01$& &{\cellcolor[gray]{0.8}}0.02 (LP Rlx)\\
     Semi-dynamic Model & 216  & 256      &  268      &    84      &  104     &  123     &  0       & 1        &   26 &147.11 &{\cellcolor[gray]{0.8}}0.09   & &{\cellcolor[gray]{0.8}}0.09 (LP Rlx)\\
     Dynamic Model      & 216  & 256      &  268      &    74      &  105     &  128     &  0       & 0       &    13 &143.41  &{\cellcolor[gray]{0.8}}0.49  & &{\cellcolor[gray]{0.8}}0.13 (LP Rlx)\\
     Perfect Information&   76     &  222     &  286      &    58      &  105     &   132    &  0       & 0        &    0   & 125.32 \\
  \bottomrule
  \end{tabular}}
  \caption{ Stochastic Models Comparison (Delay cost ratio $c_a/c_g=2$) {\bf{WITH}} TOS}\label{WithTOSTable}
\end{table*}


The optimization models are solved using Gurobi 8.1 on a laptop with 3.6 GHz processors and 32 GB RAM. The main results are listed in Tables \ref{Model_Complexity} to \ref{WithTOSTable}.
In this test example, integer solutions can be directly obtained from linear programming relaxation, for all six stochastic models and for both no-route and reroute cases.
\section{Conclusions}\label{Conclusion}
This preliminary result shows that the new binary stochastic programming model seems to be a better formulation compared with previous work. We are currently doing more numerical test and theoretical analysis.

\bibliographystyle{informs2014trsc} 




\end{document}